\documentclass[10pt]{article}
\usepackage{latexsym}

\begin{document}

\font\tenbf=cmbx10  
\font\sevenbf=cmbx7
\font\fivebf=cmbx5
\newfam\bffam
\def\bold{\fam\bffam\tenbf}
\textfont\bffam=\tenbf
\scriptfont\bffam=\sevenbf
\scriptscriptfont\bffam=\fivebf

\font\tenbbb=msbm10  
\font\sevenbbb=msbm7
\font\fivebbb=msbm5
\newfam\bbbfam
\def\Bbb{\fam\bbbfam\tenbbb}
\textfont\bbbfam=\tenbbb
\scriptfont\bbbfam=\sevenbbb
\scriptscriptfont\bbbfam=\fivebbb

\font\tengoth=eufm10  
\font\sevengoth=eufm7
\font\fivegoth=eufm5
\newfam\gothfam
\def\goth{\fam\gothfam\tengoth}
\textfont\gothfam=\tengoth
\scriptfont\gothfam=\sevengoth
\scriptscriptfont\gothfam=\fivegoth

\def\be{\begin{equation}}
\def\ee{\end{equation}}
\def\bea{\begin{eqnarray}}
\def\eea{\end{eqnarray}}

\def\G{\Gamma}
\def\k{\kappa}
\def\l{\lambda}
\def\e{\varepsilon}
\def\ra{\rangle}
\def\la{\langle}
\def\rb{\rbrack}
\def\lb{\lbrack}                              
\def\nb{\nabla}
\def\qed{\ifmmode \quad\Box
         \else $\quad\Box$ \fi}

\def\R#1{\ifmmode {\Bbb R}^{#1} 
         \else $ {\Bbb R}^{#1}$ \fi}
\def\C#1{\ifmmode {\Bbb C}^{#1} 
         \else $ {\Bbb C}^{#1}$ \fi}
\def\Z#1{\ifmmode {\Bbb Z}^{#1} 
         \else $ {\Bbb Z}^{#1}$ \fi}
\def\N{\ifmmode {\Bbb N} 
         \else $ {\Bbb N}$ \fi}
\def\H#1{\ifmmode {\Bbb H}^{#1} 
         \else $ {\Bbb H}^{#1}$ \fi}
\def\CP#1{\ifmmode {\Bbb P}^{#1}(\C{})
          \else ${\Bbb P}^{#1}(\C{})$ \fi}
\def\HP#1{\ifmmode {\Bbb P}^{#1}(\H{})
          \else ${\Bbb P}^{#1}(\H{})$ \fi}
\def\S#1{\ifmmode {\rm S}^{#1}
         \else ${\rm S}^{#1}$ \fi}
\def\spb{{\bold S}}

\def\SO{{\bold SO}}
\def\Sp{{\bold Sp}}
\def\so{{\goth so}}
\def\sl{{\goth sl}}
\def\CCl{{\Bbb C}{\bold l}}
\def\End{{\rm End}}

\def\id{{\rm id}}
\def\proj{{\rm proj}}

\def\addots{\mathinner{\mkern1mu\raise1pt\vbox{\kern7pt\hbox{.}}\mkern2mu
  \raise4pt\hbox{.}\mkern2mu\raise7pt\hbox{.}\mkern1mu}}

\renewcommand{\theequation}{\thesection.\arabic{equation}}

\def\proof{{\bf Proof.}\hskip2mm}
\def\definition{{\bf Definition.}\hskip2mm}
\def\remark{{\bf Remark.}\hskip2mm}
\def\example{{\bf Example.}\hskip2mm}
\def\leer{\par \vskip 2.7mm}

\newtheorem{Proposition}{Proposition}[section]
\newtheorem{Theorem}{Theorem}[section]
\newtheorem{Corollary}{Corollary}[section]
\newtheorem{Lemma}{Lemma}[section]

\def\nba{\tilde\nb}
\def\Oa{\tilde\Omega}
\def\cJ{{\cal J}}
\def\Da{\tilde D}
\def\Ja{\tilde J}
\def\cD{{\cal D}}

\parindent=0pt

\title{A universal lower bound for the first eigenvalue of the Dirac operator
on quaternionic K\"ahler manifolds}
\author{Wolfram Kramer\thanks{The author was
supported by the SFB 256 'Nichtlineare partielle Differentialgleichungen'}\\
Mathematisches Institut der Universit\"at Bonn\\
E-mail: {\tt kramer@rhein.iam.uni-bonn.de}}
\date{22.~August 1996}
\maketitle

\begin{abstract}
A universal lower bound for the first positive eigenvalue of the Dirac
operator on a compact quaternionic K\"ahler manifold $M$ of positive scalar
curvature is calculated. It is shown that it is equal to the first positive
eigenvalue on the quaternionic projective space. 
For this, the horizontal tangent bundle on the canonical
$\SO(3)$-bundle over $M$ is equipped with a hyperk\"ahlerian structure and the
corresponding splitting of the horizontal spinor bundle is considered. The
desired estimate is obtained by looking at hyperk\"ahlerian twistor operators
on horizontal spinors.
\end{abstract}

\section{Introduction}

On a compact Riemannian spin manifold of nonnegative scalar curvature
the first positive eigenvalue of the Dirac operator satisfies
\be
  \l^2 \geq \frac{n}{n-1}\,\frac{\inf_M \k}{4}\,,
\ee
where $\k$ denotes the scalar curvature of $M$. This result was found by 
Friedrich \cite{Friedrich1}, who has also shown that equality is attained 
if and
only if there exist at least one nontrivial Killing spinor on $M$. Killing
spinors are automatically eigenspinors for $D$ of smallest possible 
eigenvalue.

Afterwards, Hijazi \cite{Hijazi1} could show that this estimate cannot be
sharp if the manifold is K\"ahler, i.e. the lower bound cannot be attained
as eigenvalue of $D^2$. An improvement of the estimate in this case was 
done by Kirchberg \cite{Kirchberg1}, \cite{Kirchberg2}. He could show that
on a compact K\"ahler manifold $M^{2m}$ of nonnegative scalar curvature one
gets as estimates
\bea
  \l^2 & \geq & \frac{m+1}{m}\,\frac{\inf_M\k}{4} \quad \hbox{for $m$ odd}\,,
    \nonumber\\
  \l^2 & \geq & \frac{m}{m-1}\,\frac{\inf_M \k}{4} \quad \hbox{for $m$ even}\,.
\eea
In addition, Kirchberg introduced the notion of K\"ahlerian Killing spinors, 
and he could show that exactly for these spinors equality is attained.

Moreover, Hijazi's theorem \cite{Hijazi1} says that on a Riemannian spin 
manifold
$M^n$ there cannot be any nontrivial Killing spinor if there is a parallel 
$k$-form,
$k\neq 0,n$ on $M$. Hence, Friedrich's estimate is not sharp on quaternionic
K\"ahler manifolds too, because in this case there is a canonical parallel 
4-form, 
namely the Kraines form. Several attempts were made by Hijazi and
Milhorat \cite{HijaziMilhorat3}, \cite{HijaziMilhorat2} to improve the 
estimate
of the first eigenvalue also in this case. As in the case of K\"ahler 
manifolds,
they tried to introduce a suitable notion of twistor spinors and to use the 
fact that
the spinor bundle $ \spb(M)$ splits into eigenbundles $\spb_r(M)$ under 
Clifford
multiplication with the Kraines form. But up to now all estimates of the first
eigenvalue of $D^2$ depend on the concerning eigenbundle $\spb_r(M)$ in which
the spinor lives. If $M^{4m}$ is a compact spin quaternionic K\"ahler 
manifold of
positive scalar curvature, Hijazi and Milhorat could show in 
\cite{HijaziMilhorat2} that for the 
eigenvalue
$\l^2$ of an eigenspinor which lives in the bundle $\spb_0(M)$ or in the
bundle $\spb_r(M)$ for $\lb \frac{m}{2}\rb +1 \leq r \leq m $ one has the 
estimate
$\l^2 \geq \frac{m+3}{m+2}\,\frac{\k}{4}$ (there is no infimum because on
quaternionic K\"ahler manifolds the scalar curvature is always constant). 
But in all other eigenbundles the
estimate is weaker. In spite of this it was conjectured that the 
right-hand side of the estimate above
gives a universal lower bound, because it is
exactly the first eigenvalue of the spectrum of $D^2$ on the quternionic 
projective
space \cite{Milhorat1}. There was the similar situation in the Riemannian case
(Friedrich's estimate): lower bound $=$ first eigenvalue on the standard 
sphere; and
in the K\"ahler case for odd complex dimensions (Kirchberg's estimate):
lower bound $=$ first eigenvalue on the complex projective space. As a
corollary of the result in \cite{HijaziMilhorat2} it is seen that as
least for $m=2$ and $m=3$ the mentioned conjecture is correct.

In this paper it will be shown that the conjecture is true for all
quaternionic dimensions $m$:
\begin{Theorem}
Let $M^{4m}$ be a compact spin quaternionic K\"ahler manifold of positive 
scalar
curvature $\k$. Then every eigenvalue $\l$ of the Dirac operator satisfies 
the estimate
\be
  \l^2 \geq \frac{m+3}{m+2}\,\frac{\k}{4}\,.
\ee
\end{Theorem}

I would like to thank Prof.~W.~Ballmann and Prof.~C.~B\"ar for support
and encouragement and Uwe Semmelmann for many helpful discussions and
carefully reading the manuscript.

\section{Quaternionic K\"ahler manifolds, the bundle $P$}

Quaternionic K\"ahler manifolds are defined as manifolds having holonomy 
$\Sp(1)\cdot\Sp(m) = (\Sp(1) \times\Sp(m))/\Z{}_2$. Therefore there locally 
exist
three almost complex structures $J_a$, $a=1,2,3$ which satisfy the 
multiplicative
rules of imaginary unit quaternions:
\be
  J_aJ_b=\e_{abc} J_c - \delta_{ab}\,.
\ee
In addition, there are corresponding 2-forms (considered as elements in the
Clifford algebra)
\be
  \Omega_a = \frac{1}{2}\sum_{i=1}^{4m} e_i J_a e_i\,, \quad a=1,2,3\,,
\ee
which are also defined only locally. The local almost complex structures span
a three dimensional subbundle $E$ of the endomorphism bundle $\End(TM)$ of 
$TM$,
which is closed under the Levi-Civita connection. All $J_a$, $a=1,2,3$ and the
corresponding 2-forms are not parallel, but on quaternionic K\"ahler manifolds
there is a canonical parallel 4-form, the Kraines form \cite{Bonan1},
\cite{Kraines1}:
\be
  \Omega = \sum_a \Omega_a \wedge \Omega_a = \sum_a \Omega_a 
  \cdot \Omega_a + 6m\,.
\ee
Here, for the second equality, the canonical vector bundle isomorphism between
$\CCl(M)$ and $\Lambda^\ast_{\C{}}(M)$ was used and Proposition 2.3 of 
\cite{HijaziMilhorat3} was applied.
A special choice of three local almost complex structures can be identified 
with
a local frame in $E$. On the space $P$ of all of these frames the group 
$\SO(3)$
operates in a natural manner; in this way $P$ becomes a 
principal-$\SO(3)$-bundle, and let $\pi:P \to M$ denote the canonical 
projection.
The covariant derivative on $E$ induced by the Levi-Civita connection 
characterizes
a connection form $\omega$ on $P$.

To construct a metric on $P$, one considers the vector fields on $P$ 
induced by the
action of $\SO(3)$. Let $H_a$, $a=1,2,3$ be a base of $\so(3)$, orthogonal with
respect to the Killing form, which satisfies the following commutator relation:
\be
  \lb H_a,H_b \rb = 2\e_{abc} H_c\,.
\ee
The corresponding vector fields $\xi_a$, $a=1,2,3$ can be written as
\be
  \xi_a(p) = \frac{d}{dt}\Big|_{t=0}\left(p \exp(t H_a)\right)\,,
\quad p \in P\,.
\ee
Let $\omega_a$ be the 1-form dual to $\xi_a$ which annihilates all
vectors that are horizontal w.r.t. $\omega$, i.e. $\omega$ can be
written as
\be
  \omega = \sum_a c_a \omega_a
\ee
for some $c_a \in \R{}$. The metric on $P$ can now be defined by
\be
  g_P:= \pi^{\ast} g_M + \sum_a \omega_a\otimes \omega_a\,.
\ee
Hence, there is a orthogonal splitting of the tangent bundle of $P$ into a
horizontal and a vertical part:
\be
  TP = T_HP \oplus T_VP\,.
\ee
In the future only the horizontal bundle will be of special interest. As 
abbreviation
$HP = T_HP$ will be used. The idea of the proof is to lift all calculations 
from
the quaternionic K\"ahler manifold to the bundle $P$ as exactly as possible 
but with
the difference that on $P$ there are now globally defined complex structures 
which
in addition are parallel w.r.t. the horizontal connection on $HP$. That means
that the horizontal bundle will be equipped with a hyperk\"ahler structure.

At this point it is convenient to introduce some conventions.
If $X \in TM$ is a vector on $M$,
let $X^\ast \in HP$ denote its horizontal lifting to $P$. Explicitly one has
$X^\ast \perp T_VP$. The horizontal connection $\nba$ on $HP$ is defined by
\be
  \nba_{X^\ast}Y^\ast := (\nb_XY)^\ast \qquad \hbox{and} \qquad
  \nba_VY^\ast = 0
\ee
for all $V \in T_VP$.

Now one considers the pull-back $\pi^\ast \spb(M)$ of the spinor bundle 
$\spb(M)$
onto $P$, which is isomorphic to the spinor bundle $\spb(HP)$ associated to 
$HP$.
It is a $\CCl(HP)$-module in a natural manner, where the Clifford 
multiplication
operates by
\be
  X^\ast \pi^\ast\psi(p) := \pi^\ast(X\psi)(p), \qquad \psi \in \G(\spb(M))\,.
\ee
Moreover, $\spb(HP)$ is equipped with the pull-back connection on $M$:
\be
  \nba_{X^\ast}\pi^\ast\psi(p) = \pi^\ast(\nb_X\psi)(p)\,.
\ee
On sections of $\spb(HP)$ a horizontal Dirac operator $\Da$ is defined by
\leer 
\definition
\be
  \Da\psi = \sum_{i=1}^{4m} e_i^\ast \nba_{e_i^\ast} \psi, \qquad \psi\in 
\G(\spb(HP))\,.
\ee
Here, $\{e^\ast_i\}$ is a orthonormal base of $HP$.
\leer

The key point of the following investigations is
\begin{Proposition} \label{wennDdannDa}
Let $\psi \in \G(\spb(M))$ be an eigenspinor of $D$ with eigenvalue $\l$. Then
$\pi^\ast\psi \in \G(\spb(HP))$ is an eigenspinor of $\Da$ with the same
eigenvalue $\l$.
\end{Proposition}
\proof
The proof is straightforward because of the definition of the horizontal
connection $\nba$. \qed
\leer

On $HP$ there are naturally defined global almost complex structures. A
point $p\in P$ is a frame of three almost complex structures over a point
$m\in M$:
\be
  p = (J_1(p),J_2(p),J_3(p))\,.
\ee
\leer
\definition
Let $J_1,J_2,J_3$ be the three almost complex structures on $HP$, defined by
\be
  J_a X^\ast(p) := (J_a(p)X)^\ast(p)\,.
\ee

\begin{Proposition}
$J_1$, $J_2$ and $J_3$ on $P$ are parallel w.r.t. $\nba$.
\end{Proposition}

\proof
It has to be shown that $\nba_{X^\ast}J_1 = 0 $ for all $X^\ast\in HP$.
But this is seen at once, because by definition the connection form on
$P$ satisfies $\omega(X^\ast)=0$. \qed
\leer

Therefore, three K\"ahler forms can be defined on $P$ which are denoted in
the following by $\Oa_a$, $a=1,2,3$:
\be
  \Oa_a = \frac{1}{2}\sum_{i=1}^{4m} e_i J_a e_i\,,
\ee
where $\{e^\ast_i\}$ is a orthonormal base of $HP$. In the same manner 
there exists a horizontal Kraines form on $P$:
\be
  \Oa := \sum_a \Oa_a \wedge \Oa_a = \sum_a \Oa_a \cdot \Oa_a + 6m\,.
\ee
One has $ \Oa = \pi^\ast \Omega$, and $\Oa$ is parallel w.r.t. $\nba$. 
Hence $HP$ is equipped with a hyperk\"ahler
structure w.r.t. $\nba$.

\section{Splitting of the horizontal spinor bundle}
\label{Zerlegung}

Since $HP$ is the only bundle which is dealt with, the star 
$\hbox to 0pt{}^\ast$
that denotes horizontal liftings of vectors will be omitted, and the 
short notation
$\spb := \spb(HP)$ will be used.

The spinor bundle $\spb(M)$ on $M$ splits into eigenbundles of the 
Kraines form;
this was shown by Hijazi and Milhorat \cite{HijaziMilhorat1}. 
This splitting is carried over to $\spb$
at once. Besides of the horizontal Krainesform $\Oa$ on $P$ one can choose 
one of
the three horizontal K\"ahlerforms which in the following will always be 
denoted by
$\Oa_1$. $\Oa$ and $\Oa_1$ are parallel w.r.t. $\nba$. Because of
$\lb \Oa,\Oa_1\rb =0$ there is in addition to the mentioned splitting of $\spb$
a  decomposition into eigenbundles of $\Oa_1$. For further investigation it is
necessary to look at representations of $\sl(2,\C{})$.

The horizontal K\"ahler forms $\Oa_a$ satisfy the following commutator 
relations:
\be
  \lb \Oa_a,\Oa_b \rb = 4 \e_{abc}\Oa_c\,.
\ee
One considers in $\CCl(HP)$ the forms
\be
  O_a := \frac{i}{2} \Oa_a
\ee
and
\be
  O_1^+:= \frac{1}{2}(O_2 + i O_3)\,, \quad 
  O_1^-:= \frac{1}{2}(O_2 - i O_3)\,.
\ee
It can be verified at once that this is a representation of $\sl(2,\C{})$:
\bea
  \lb O_1,O_1^+ \rb & = & 2 O_1^+ \nonumber\\
  \lb O_1,O_1^- \rb & = & -2 O_1^- \nonumber\\
  \lb O_1^+,O_1^- \rb & = & O_1\,.
\eea
The corresponding Casimir operator is easily written down. W.r.t. the
Killing form, $\frac{1}{8}O_1$ is dual to $O_1$, and $\frac{1}{4}O_1^+$ resp.
$\frac{1}{4}O_1^-$ is dual to $O_1^-$ resp. $O_1^+$. Hence the Casimir operator
is given by
\bea 
  C &=& \frac{1}{8}O_1 O_1 + \frac{1}{4}O_1^+ O_1^- + \frac{1}{4}O_1^- O_1^+
   = \frac{1}{8} \Big( \sum_a O_a O_a \Big) \nonumber\\
  &=& -\frac{1}{32}(\Oa -6m)\,.
\eea
In general, the Casimir operator operates on an irreducible representation of 
highest weight $\mu$
by multiplication with $\|\mu + \rho \|^2 - \|\rho\|^2$,
where $\rho$ denotes the half sum of positive roots. The scalar product on
the space of weights is defined by $\la \mu,\nu \ra := B(t_\mu,t_\nu)$,
where $t_\mu$ is the uniquely determined element of the Cartan subalgebra 
${\goth h}$ with $ B(t_\mu,h) = \mu(h)$ for all $h\in {\goth h}$. In the 
special
case $\sl(2,\C{})$ (the Cartan subalgebra has dimension 1 and is spanned 
by $O_1$)
$\mu$ is simply a natural number $r$. Let $V$ be a irreducible representation
of $\sl(2,\C{})$ with highest weight $r$, then the vector of highest weight
satisfies
\be
  O_1 v_r = r v_r\,,
\ee
and $C$ operates on $V$ by multiplication with $\frac{1}{8}r(r+2)$.

Using this the splitting $\spb(HP)$ into eigenbundles w.r.t. the operation 
of $\Oa$ and $\Oa_1$ is determined:
\bea
  \spb &=& \bigoplus_{r=0}^m \spb_r, \qquad \Oa|_{\spb_r} = 
  (6m - 4r(r+2))\id\,,\nonumber\\
  \spb_r &=& \bigoplus_{s=0}^r \spb_r^{(s)}, \qquad \Oa_1|_{\spb_r^{(s)}} 
  = i(2r-4s) \id\,.
\eea
To summarize, one gets the following facts: if the partial bundle $\spb_r$ with
eigenvalue $6m-4r(r+2)$ w.r.t. $\Oa$ is considered, there is an additional
splitting into smaller bundles, which are eigenbundles of the Clifford 
multiplication
with $\Oa_1$. Between the eigenvalues of $\Oa_1$ on these partial bundles there
are gaps of absolute value 4.

On the other side it is known that like spinor bundles on K\"ahler manifolds
the whole of $\spb$ splits into eigenbundles of $\Oa_1$ with eigenvalues 
$i(2m-2k)$:
$\spb = \bigoplus_{k=0}^{2m} \spb^k$, $k = 0,\dots,2m$. Explicitly one has
\be
  \spb_r^{(s)} = \spb_r \cap \spb^{m-r+2s}\,.
\ee
To avoid inconveniences with notations the definition
\be
  \spb_r^k := \spb_r^{(\frac{k+r-m}{2})} = \spb_r \cap \spb^k
\ee
for ${\frac{k+r-m}{2}} \in \N{}_0$ will be used. In the future the 
appearence of
$\spb_r^k$ should be interpreted in the sense that all which is said should be
ignored if $k$ does not satisfy the integrability condition or if the bundle
$\spb_r^k$ does not exist at all (if e.g. $k<0$).

All that was said can be clarified by a picture:

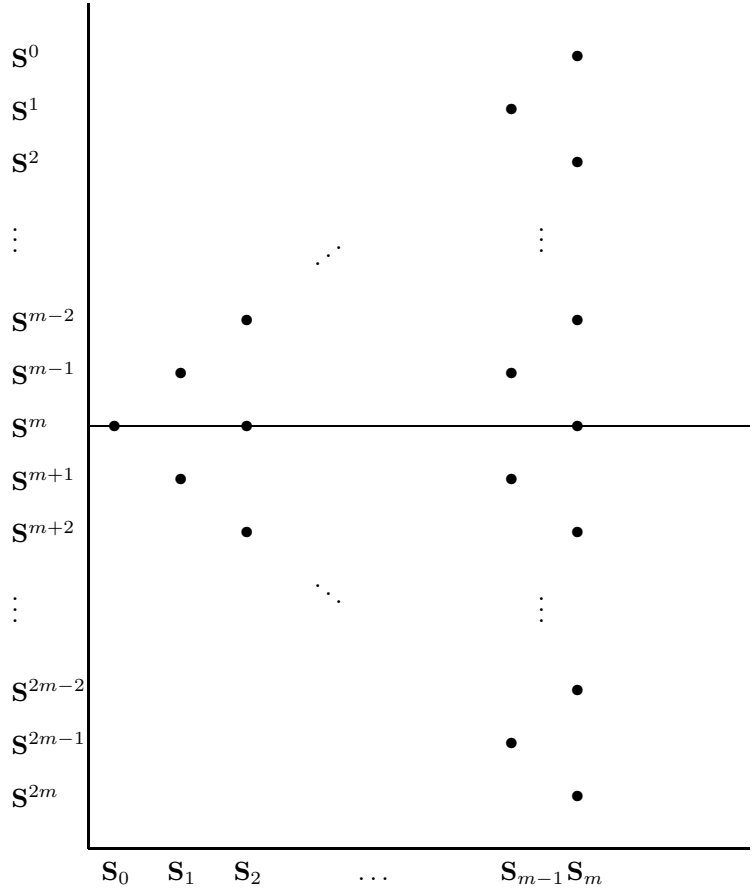
\begin{center}
\begin{picture}(250,350)(-50,-185)

\put(0,0){\circle*{4}}
\put(25,20){\circle*{4}}
\put(25,-20){\circle*{4}}
\put(50,40){\circle*{4}}
\put(50,0){\circle*{4}}
\put(50,-40){\circle*{4}}

\put(75,60){$\addots$}
\put(75,-68){$\ddots$}

\put(160,66){\vdots}
\put(160,-74){\vdots}
\put(150,120){\circle*{4}}
\put(150,20){\circle*{4}}
\put(150,-20){\circle*{4}}
\put(150,-120){\circle*{4}}
\put(175,140){\circle*{4}}
\put(175,100){\circle*{4}}
\put(175,0){\circle*{4}}
\put(175,40){\circle*{4}}
\put(175,-40){\circle*{4}}
\put(175,-100){\circle*{4}}
\put(175,-140){\circle*{4}}

\put(-10,0){\line(1,0){250}}
\put(-10,-160){\line(0,1){320}}
\put(-10,-160){\line(1,0){250}}

\put(-5,-172){$\spb_0$}
\put(20,-172){$\spb_1$}
\put(45,-172){$\spb_2$}
\put(92,-172){\dots}
\put(146,-172){$\spb_{m-1}$}
\put(171,-172){$\spb_m$}

\put(-39,136){$\spb^0$}
\put(-39,116){$\spb^1$}
\put(-39,96){$\spb^2$}
\put(-39,66){\vdots}
\put(-39,36){$\spb^{m-2}$}
\put(-39,16){$\spb^{m-1}$}
\put(-39,-4){$\spb^m$}
\put(-39,-24){$\spb^{m+1}$}
\put(-39,-44){$\spb^{m+2}$}
\put(-39,-74){\vdots}
\put(-39,-104){$\spb^{2m-2}$}
\put(-39,-124){$\spb^{2m-1}$}
\put(-39,-144){$\spb^{2m}$}

\end{picture}
\end{center}

The Clifford multiplication with vectors in $HP$ is not compatible with the
splitting. In fact, a spinor lying completely in one distinguished eigenbundle
of $\Oa$ resp. $\Oa_1$ will be carried into the direct sum of the two 
neighbouring
eigenbundles:

\begin{Lemma}
\bea
  \mu: HP \otimes \spb_r & \to & \spb_{r-1}\oplus \spb_{r+1} \quad \hbox{resp.}
    \nonumber\\
  \mu: HP \otimes \spb^k & \to & \spb^{k-1}\oplus \spb^{k+1}\,.
\eea
\end{Lemma}

\proof
See \cite{Kirchberg2}, \cite{HijaziMilhorat2}. \qed
\leer

Summarizing these facts one gets
\be \label{aufspmu}
  \mu:= \mu^{++}\oplus\mu^{+-}\oplus\mu^{-+}\oplus\mu^{--}:
 HP \otimes \spb^k_r \to \spb^{k+1}_{r+1}\oplus \spb^{k-1}_{r+1}
   \oplus \spb^{k+1}_{r-1}\oplus \spb^{k-1}_{r-1}
\ee
(here every non-existing summand in the Whitney sum is to be omitted). To
project Clifford multiplication on every summand, one has to follow the ideas
of Kirchberg and Hijazi and to combine them. The following lemmata are given
without proofs. For the proof of Lemma \ref{lemma1} see \cite{Kirchberg2} and 
for the proofs of Lemmata \ref{lemma2} and \ref{lemma3} see
\cite{HijaziMilhorat2}.
If the bundle $HP$ is complexified and the Clifford multiplication is linearly
extended, the following facts are easily proven:
\begin{Lemma}\label{lemma1}
\bea
  q^+(X) = \frac{1}{2}(X +i J_1 X) &:& \spb^k \to \spb^{k+1} \nonumber\\
  q^-(X) = \frac{1}{2}(X -i J_1 X) &:& \spb^k \to \spb^{k-1}\,. \qed
\eea
\end{Lemma}
Hijazi \cite{HijaziMilhorat2} invented the operator
$\cJ$, which is defined by $ \cJ(X) := \sum_a \Oa_aJ_a(X) +3X $
for the treatment of the splitting w.r.t. $\Oa$:
\begin{Lemma}\label{lemma2}
\bea
  \lb \Oa, X \rb &=& 4 \cJ(X) \nonumber\\
  \lb \Oa, \cJ(X)\rb &=& -8\cJ(X) + 12 X -4X(\Oa -6m)\,. \qed
\eea
\end{Lemma}
With help of this one can define similar projectors in the case of the 
splitting 
in eigenbundles w.r.t. the Kraines form:
\begin{Lemma}\label{lemma3}
\bea
  p_r^+(X) = \frac{1}{4(r+1)}\big((2r+1)X- \cJ(X)\big) &:& 
   \spb_r \to \spb_{r+1}
   \nonumber\\
  p_r^-(X) = \frac{1}{4(r+1)}\big((2r+3)X+ \cJ(X)\big) &:& 
  \spb_r \to \spb_{r-1}\,.
\qed
\eea
\end{Lemma}

Let $\{e_j\}$ be an orthogonal base of $HP$ with the property 
$e_{2j}= J_1 e_{2j-1}$,
$j=1,\dots,2m$. $f_j = q^-(e_{2j-1})$ and $\bar f_j = q^+(e_{2j-1})$ form the
corresponding complex base. The operation of the projectors can be interpreted
in a nice graphical way:

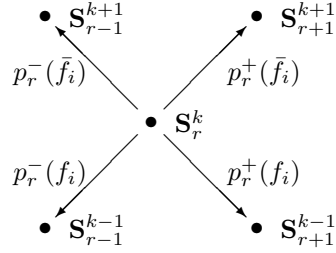
\begin{center}
\begin{picture}(120,120)(-60,-60)
\put(0,0){\circle*{4}}
\put(9,-4){$\spb_r^k$}
\put(40,40){\circle*{4}}
\put(49,36){$\spb_{r+1}^{k+1}$}
\put(5,5){\vector(1,1){31}}
\put(29,16){$p_r^+(\bar f_i)$}
\put(40,-40){\circle*{4}}
\put(49,-44){$\spb_{r+1}^{k-1}$}
\put(5,-5){\vector(1,-1){31}}
\put(29,-22){$p_r^+(f_i)$}
\put(-40,40){\circle*{4}}
\put(-31,36){$\spb_{r-1}^{k+1}$}
\put(-5,5){\vector(-1,1){31}}
\put(-52,16){$p_r^-(\bar f_i)$}
\put(-40,-40){\circle*{4}}
\put(-31,-44){$\spb_{r-1}^{k-1}$}
\put(-5,-5){\vector(-1,-1){31}}
\put(-52,-22){$p_r^-(f_i)$}
\end{picture}
\end{center}

\section{Hyperk\"ahlerian Twistor operators}

The twistor operator ${\cal D}$ on a Riemannian spin manifold $M^n$ is by 
definition
the composition of the covariant differential $\nb$ and followed by the 
orthogonal
projection onto the kernel of the Clifford multiplication:
\bea
  {\cal D}\psi = \proj_{\ker \mu}(\nb \psi)
  &=& \sum_i e_i\otimes \nb_{e_i}\psi + 
  \frac{1}{n}\sum_{i,j}e_i\otimes e_i e_j \nb_{e_j}\psi
 \nonumber\\
  &=& \sum_i e_i\otimes \nb_{e_i}\psi + 
  \frac{1}{n}\sum_{i,j}e_i\otimes e_i D\psi
\eea
for an arbitrary orthogonal base $\{e_i\}$.
Lower estimates for the first eigenvalue of the Dirac operator are established
by considering the inequality $\|{\cal D} \psi\|^2 \geq 0$. But if the manifold
$M$ carries additional structure, this is not sufficient, as shown by Hijazi 
\cite{Hijazi1} in the case of K\"ahler manifolds. Here it was necessary to 
split
the Clifford multiplication in a similar way as above and to define partial
twistor operators.

This approach will be used also in the case of a hyperk\"ahler structure
on $HP$. Corresponding to the splitting (\ref{aufspmu}) of Clifford 
multiplication
the horizontal Dirac operator can be splitted into four parts:
\be
  \Da = D^{++} + D^{+-} + D^{-+} + D^{--}\,,
\ee
which, restricted to $\G(\spb_r^k)$, have the following form:
\bea
  D^{++} &=& 2 \sum_{j=1}^{2m} p_r^+(\bar f_j) \nba_{f_j} \nonumber\\
  D^{+-} &=& 2 \sum_{j=1}^{2m} p_r^+(f_j) \nba_{\bar f_j} \nonumber\\
  D^{-+} &=& 2 \sum_{j=1}^{2m} p_r^-(\bar f_j) \nba_{f_j} \nonumber\\
  D^{--} &=& 2 \sum_{j=1}^{2m} p_r^-(f_j) \nba_{\bar f_j}\,.
\eea
Therefore the partial twistor operators are easily written down
(again as above, restricted to $\G(\spb_r^k)$):
\begin{Lemma}
\bea
  {\cal D}^{++} &=& \sum_{j=1}^{2m}\Big(
  p_r^+(\bar f_j) \otimes \nba_{f_j} 
  \nonumber\\
  && \phantom{\sum_{j=1}^{2m}\Big(} - \frac{1}{A_{r+1,k+1}^{++}}
    p_r^+ (\bar f_j)\otimes p_{r+1}^-(f_j)D^{++}
    \Big)\nonumber\\
  {\cal D}^{+-} &=& \sum_{j=1}^{2m}\Big(
  p_r^+(f_j) \otimes \nba_{\bar f_j} 
  \nonumber\\
  && \phantom{\sum_{j=1}^{2m}\Big(} - \frac{1}{A_{r+1,k-1}^{+-}}
    p_r^+(f_j)\otimes p_{r+1}^-(\bar f_j)D^{+-}
    \Big)\nonumber\\
  {\cal D}^{-+} &=& \sum_{j=1}^{2m}\Big(
  p_r^-(\bar f_j) \otimes \nba_{f_j} 
  \nonumber\\
  && \phantom{\sum_{j=1}^{2m}\Big(} - \frac{1}{A_{r-1,k+1}^{-+}}
    p_r^-(\bar f_j)\otimes p_{r-1}^+(f_j)D^{-+}
    \Big)\nonumber\\
  {\cal D}^{--} &=& \sum_{j=1}^{2m}\Big(
  p_r^-(f_j) \otimes \nba_{\bar f_j} 
  \nonumber\\
  && \phantom{\sum_{j=1}^{2m}\Big(} - \frac{1}{A_{r-1,k-1}^{--}}
    p_r^-(f_j)\otimes p_{r-1}^+(\bar f_j)D^{--}\Big)
\eea
with
\bea
   A_{r,k}^{\pm+} &=& \sum_{j=1}^{2m}p_{r\mp 1}^\pm(\bar f_j)p_{r}^\mp(f_j)
    \big|_{\spb_r^k}
 \qquad\hbox{and} \nonumber\\
   A_{r,k}^{\pm-} &=& \sum_{j=1}^{2m}p_{r\mp 1}^\pm(f_j)p_{r}^\mp(\bar f_j)
    \big|_{\spb_r^k}\,.
\eea
\end{Lemma}
\proof
It is easily verified that the twistor operators defined above are lying in
the kernel of the Clifford multiplication. It remains to prove the 
orthogonality
of the projection. But this is clear by the observation that e.g. 
$p_r^+(\bar f_j)$
is the adjoint of $p_{r+1}^-(\bar f_j)$ w.r.t. the fibrewise scalar
product on spinors. \qed
\leer

In the next step the absolute value of the twistor operators, applied 
to a spinor,
will be calculated. It is therefore sufficient to consider the following
reduced twistor operators:
\bea
  {\cal D}_X^{\pm+} &=& 
    \nba_{q^-(X)} - \sum_{k=1}^{2m}\frac{1}{A_{r\pm 1,k+1}^{\pm+}}
    p_{r\pm 1}^\mp(q^-(X))p_r^\pm(\bar f_k) \nba_{f_k}
    \nonumber\\
  {\cal D}_X^{\pm-} &=& 
    \nba_{q^+(X)} - \sum_{k=1}^{2m}\frac{1}{A_{r\pm 1,k-1}^{\pm-}}
    p_{r\pm 1}^\mp(q^+(X))p_r^\pm(f_k) \nba_{\bar f_k}
    \,.
\eea

In the rest of this section the constants $A_{r,k}^{\pm\pm}$ are calculated
explicitly. In order to do this, some technical lemmata are necessary. In the
following calculations dealing with the complex structures $J_a$, one often 
has to 
distinguish the cases $a=1$ and $a\neq 1$. Therefore the convention will 
be used
that summation over $a'$ means summation over $a=2$ and $a=3$ but not over 
$a=1$.
Moreover, w.r.t. the summation over the indices $a$ and $b$ the Einstein
convention is used to avoid too many sums.

\begin{Lemma}
\be
  \sum_j J_{a'}f_j J_{a'}\bar f_j = \sum_j \bar f_j f_j\,.
\ee
\end{Lemma}
\proof
\be
  \sum_{j} J_{a'}f_j J_{a'}\bar f_j = \frac{1}{4}\sum_{j}
  (J_{a'}e_{2j-1}-iJ_{a'}e_{2j})(J_{a'}e_{2j-1}+iJ_{a'}e_{2j})\,.
\ee
A new base is defined by $e_j' = (-1)^{j+1} J_{a'}e_j$. Hence,
\be
  J_1 e'_{2j-1} = J_1J_{a'}e_{2i-1}=-J_{a'}J_1e_{2j-1} 
  = -J_{a'}e_{2j}= e'_{2j}\,.
\ee
Therefore
\be
  \sum_{j} J_{a'}f_j J_{a'}\bar f_j = \frac{1}{4}\sum_{j}
  (e'_{2j-1}+ie'_{2j})(e'_{2j-1}-ie'_{2j})
  =\frac{1}{4}\sum_j \bar f'_j f'_j\,. \qed
\ee
\leer

\begin{Lemma}
\be
  J_{a'}f_j\bar f_j = - \bar f_j J_{a'}f_j\,.
\ee
\end{Lemma}
\proof
\be
  (J_{a'}e_{2j-1}-iJ_{a'}e_{2j})(e_{2j-1}+ie_{2j}) 
  = - (e_{2j-1}+ie_{2j})(J_{a'}e_{2j-1}-iJ_{a'}e_{2j})\,,
\ee
since $J_{a'}e_{2j-1}$ and $J_{a'}e_{2j}$ are orthogonal to $e_{2j-1}$ and
$e_{2j}$. \qed
\leer

\begin{Lemma}
\bea
  \sum_j f_jJ_{a'}\bar f_j &=& \frac{1}{2}\Oa_{a'} + 
  \frac{i}{2}\e_{a'1b'}\Oa_{b'} 
    \nonumber\\
  \sum_j \bar f_jJ_{a'} f_j &=& \frac{1}{2}\Oa_{a'} - 
  \frac{i}{2}\e_{a'1b'}\Oa_{b'}\,.
\eea
\end{Lemma}
\proof
\bea
  \sum_j f_jJ_{a'}\bar f_j &=& 
  \frac{1}{4}\sum_j(e_{2j-1}-ie_{2j})J_{a'}(e_{2j-1}+ie_{2j})
   \nonumber\\
   &=& \frac{1}{4}\sum_j (e_{2j-1}J_{a'}e_{2j-1} + ie_{2j}J_{a'}e_{2j-1}
      \nonumber\\
    &&\phantom{\frac{1}{4}\sum_j (}
    -i e_{2j-1}J_{a'}e_{2j} + e_{2j}J_{a'}e_{2j}) \nonumber\\
   &=& \frac{1}{2}\Oa_{a'} + \frac{i}{4}\sum_j (e_{2j-1}J_{a'}e_{2j} 
      - e_{2j}J_{a'}e_{2j-1})  \nonumber\\
   &=& \frac{1}{2}\Oa_{a'} + \frac{i}{4}\sum_j (e_{2j-1}J_{a'}J_1e_{2j-1} 
      + e_{2j}J_{a'}J_1e_{2j})  \nonumber\\
   &=& \frac{1}{2}\Oa_{a'} + \frac{i}{2}\e_{a'1b'}\Oa_{b'}\,.
\eea
The second equation is proven analogously. \qed

\begin{Lemma}
\bea
  {\cal J}(f_j) &=& \Oa_{a'}J_{a'} f_j + (3 + i\Oa_1)f_j
     \nonumber\\
  {\cal J}(\bar f_j) &=& \Oa_{a'}J_{a'} \bar f_j + (3 - i\Oa_1)\bar f_j\,.
\eea
\end{Lemma}
\proof
\bea
  {\cal J}(f_j) &=& \Oa_aJ_af_j + 3f_j = \Oa_{a'}J_{a'}f_j + 
  \Oa_1J_1 f_j+ 3f_j 
     \nonumber\\
  &=& \Oa_{a'}J_{a'}f_j + (3+ i\Oa_1)f_j
\eea
and similarly
\bea
  {\cal J}(\bar f_j) &=& \Oa_aJ_a \bar f_j + 3\bar f_j = 
    \Oa_{a'}J_{a'}\bar f_j + \Oa_1J_1 \bar f_j+ 3\bar f_j 
     \nonumber\\
  &=& \Oa_{a'}J_{a'}\bar f_j + (3- i\Oa_1)\bar f_j\,. \qed
\eea
\leer
In order to make the following expressions more being able to be handled, the
notations
\bea
  L &=& \sum_{j}\Oa_{a'}f_jJ_{a'}\bar f_j \nonumber\\
  \bar L &=& \sum_{j}\Oa_{a'}\bar f_jJ_{a'}f_j
\eea
are introduced.

\begin{Lemma}
\bea
  {\cal J}(f_j)\bar f_j &=& -\bar L + (3+ i \Oa_1)f_j\bar f_j \nonumber\\
  {\cal J}(\bar f_j)f_j &=& - L + (3- i \Oa_1)\bar f_j f_j \nonumber\\
  f_j{\cal J}(\bar f_j) &=& L + (1- i \Oa_1) f_j \bar f_j -4\bar f_j f_j 
  \nonumber\\
  {\cal J}(\bar f_j)f_j &=& \bar L + (1+ i \Oa_1)\bar f_j f_j 
   - 4f_j \bar f_j\,.
\eea
\end{Lemma}
\proof
For example:
\bea
  {\cal J}(f_j)\bar f_j &=& (\Oa_{a'}J_{a'}f_j + 
  (3+ i \Oa_1)f_j)\bar f_j \nonumber\\
   &=& -\Oa_{a'}\bar f_jJ_{a'}f_j + (3+ i \Oa_1)f_j\bar f_j \nonumber\\
   &=& - \bar L + (3+ i \Oa_1)f_j\bar f_j \\
   f_j{\cal J}(\bar f_j) 
   &=& f_j(\Oa_{a'}J_{a'}\bar f_j + (3- i \Oa_1)\bar f_j)\nonumber\\
   &=& \Oa_{a'}f_j J_{a'}\bar f_j - 2J_{a'}f_jJ_{a'}\bar f_j 
  + (3- i \Oa_1)f_j\bar f_j
     + 2i J_1 f_j\bar f_j \nonumber\\
  &=& L + (1- i \Oa_1) f_j \bar f_j -4\bar f_j f_j\,. \qed
\eea
\leer

\begin{Lemma}
\bea
  \sum_j{\cal J}(f_j){\cal J}(\bar f_j) &=&
  \sum_j\Big(-12\bar f_jf_j + \Oa_{a'}\Oa_{a'}
  \bar f_jf_j  + (-1+i\Oa_1)L 
   - (1-i\Oa_1)\bar L 
   \nonumber\\ && \phantom{\sum_j \Big(}
  +4L +(3+i\Oa_1)(1-i\Oa_1)
  f_j\bar f_j \Big)\,.
\eea
\end{Lemma}
\proof
\bea
  \sum_j{\cal J}(f_j){\cal J}(\bar f_j) &=& \sum_j
  \Big( \Oa_{a'}J_{a'}f_j + (3+ i\Oa_1)f_j \Big)
  \Big( \Oa_{a'}J_{a'}\bar f_j + (3- i\Oa_1)\bar f_j \Big)
   \nonumber\\
  &=& \sum_j \Big(\Oa_{a'}J_{a'}f_j\Oa_{b'}J_{b'}\bar f_j
     + \Oa_{a'}J_{a'}f_j(3- i\Oa_1)\bar f_j
     \nonumber\\&& \phantom{\sum_j \Big(}
     + (3+ i\Oa_1)f_j\Oa_{a'}J_{a'}\bar f_j 
     + (3+ i\Oa_1)f_j(3- i\Oa_1)\bar f_j\Big)
     \nonumber\\
  &=& \sum_j \Big( \Oa_{a'}\Oa_{b'} J_{a'}f_jJ_{b'}\bar f_j
     - 2\Oa_{a'}J_{b'}J_{a'}f_jJ_{b'}\bar f_j
     \nonumber\\&& \phantom{\sum_j \Big(}
     + \Oa_{a'}(3- i\Oa_1)J_{a'}f_j\bar f_j
     + 2i \Oa_{a'}J_1J_{a'}f_j\bar f_j
     \nonumber\\&& \phantom{\sum_j \Big(}
     + (3+ i\Oa_1)\Oa_{a'}f_jJ_{a'}\bar f_j
    - 2 (3+ i\Oa_1)J_{a'}f_jJ_{a'}\bar f_j
     \nonumber\\ && \phantom{\sum_j \Big(}
      + (3+ i\Oa_1)(1- i\Oa_1)f_j\bar f_j \Big)\nonumber\\
  &=& \sum_j \Big( 4i\Oa_1\bar f_j f_j 
       + \Oa_{a'}\Oa_{a'}\bar f_j f_j
       + (3- i\Oa_1)\Oa_{a'}J_{a'}f_j\bar f_j
       \nonumber\\&& \phantom{\sum_j \Big(}
       - 4i \e_{a'1b'}\Oa_{b'}J_{a'}f_j\bar f_j
            + 2 \Oa_{a'}J_{a'}f_j\bar f_j
       \nonumber\\&& \phantom{\sum_j \Big(}
       + (3+ i\Oa_1)L - 4(3+ i\Oa_1)\bar f_j f_j
        \nonumber\\ && \phantom{\sum_j \Big(}
    + (3+ i\Oa_1)(1- i\Oa_1)f_j\bar f_j \Big)\nonumber\\
  &=& \sum_j \Big(
     - 12\bar f_jf_j + (-1+i\Oa_1)L 
     + \Oa_{a'}\Oa_{a'}\bar f_j f_j 
      \nonumber\\ && \phantom{\sum_j \Big(}
    - (3- i\Oa_1)\Oa_{a'}\bar f_jJ_{a'} f_j
     + 4i \e_{a'1b'}\Oa_{b'}J_1\bar f_j J_{a'}J_1f_j
     \nonumber\\ && \phantom{\sum_j \Big(}
     - 2\bar L 
     +(3+ i\Oa_1)(1- i\Oa_1)f_j\bar f_j  \Big)\nonumber\\
  &=& \sum_j \Big(
     - 12\bar f_jf_j + \Oa_{a'}\Oa_{a'}\bar f_j f_j 
     + (-1+i\Oa_1)L + 4L 
     \nonumber\\&& \phantom{\sum_j \Big(}
     -(5-i\Oa_1)\bar L + 4\bar L
     + (3+ i\Oa_1)(1- i\Oa_1)f_j\bar f_j  \Big)\nonumber\\
  &=& \sum_j \Big(
     - 12\bar f_jf_j + \Oa_{a'}\Oa_{a'}\bar f_j f_j 
     + (-1+i\Oa_1)L - (1-i\Oa_1) \bar L 
     \nonumber\\ && \phantom{\sum_j \Big(}
     +4L +(3+ i\Oa_1)(1- i\Oa_1)f_j\bar f_j  \Big)\,.
  \qed
\eea
\leer

Analogously it can be calculated:
\begin{Lemma}
\bea
  \sum_j{\cal J}(\bar f_j){\cal J}(f_j) &=&
  \sum_j\Big(-12 f_j\bar f_j + \Oa_{a'}\Oa_{a'}
  f_j\bar f_j  -(1+i\Oa_1)L 
  + (-1-i\Oa_1)\bar L     
  \nonumber\\ && \phantom{\sum_j \Big(}
   +4\bar L +(3-i\Oa_1)(1+i\Oa_1)
    \bar f_j f_j \Big)\,. \qed
\eea
\end{Lemma}

\begin{Lemma}
After restriction to $\spb_r^k$ one gets:
\bea
  \Oa_1|_{\spb_r^k} &=& (2m-2k)\id \nonumber\\
  \Oa|_{\spb_r^k} &=& (6m- 4r(r+2))\id \nonumber\\
  L|_{\spb_r^k} &=& \big(-2r(r+2) + (m-k)(2m-2k+4)\big)\id \nonumber\\
  \bar L|_{\spb_r^k} &=& \big(-2r(r+2) + (m-k)(2m-2k-4)\big)\id\,.
\eea
\end{Lemma}
\proof
Only the last two equations are not obvious.
\bea
  L &=& \sum_{j}\Oa_{a'}f_jJ_{a'}\bar f_j \nonumber\\
    &=& \sum_{j}\Oa_{a'}\big( \frac{1}{2}\Oa_{a'}
       + \frac{i}{2}\e_{a'1b'}\Oa_{b'}\big)\nonumber\\
    &=& \frac{1}{2}\big((\Oa -6m)-\Oa_1\Oa_1) -2i \Oa_1\,.
\eea
If $\Oa$ and $\Oa_1$ are replaced by the numerical values, the desired
result is obtained. The calculation of $\bar L$ is almost the same. \qed
\leer

With these lemmata, one is able to calculate the projector sums
\bea
  \sum_{j=1}^{2m} p^-_{r+1}(f_j)p^+_r(\bar f_j) &=&
    \frac{1}{16(r+2)(r+1)}\Big((2r+5)(2r+1)f_j\bar f_j 
   \nonumber\\ &&
   + (2r+1){\cal J}(f_j)\bar f_j-(2r+5)f_j{\cal J}(\bar f_j) 
   - {\cal J}(f_j){\cal J}(\bar f_j)\Big)
   \nonumber\\ &&
\eea
(and the corresponding expressions for the other three sums) and finally the
numerical values of the constants $A_{r,k}^{\pm\pm}$. These are long and ugly
calculations so the author has used a computer program.

\begin{Proposition} \label{propA}
\bea
  A_{r,k}^{--} &=&\sum_{j=1}^{2m}p_{r+1}^-(f_i)p_r^+(\bar f_i)\big|_{\spb_r^k}
    = \frac{(-m+r)(2+k-m+r)}{2(r+1)}\nonumber\\
  A_{r,k}^{+-} &=&\sum_{j=1}^{2m}p_{r-1}^+(f_i)p_r^-(\bar f_i)\big|_{\spb_r^k}
    = \frac{(k-m-r)(2+m+r)}{2(r+1)}\nonumber\\
  A_{r,k}^{-+} &=&\sum_{j=1}^{2m}p_{r+1}^-(\bar f_i)p_r^+(f_i)\big|_{\spb_r^k}
    = \frac{(-m+r)(2-k+m+r)}{2(r+1)}\nonumber\\
  A_{r,k}^{++} &=&\sum_{j=1}^{2m}p_{r-1}^+(\bar f_i)p_r^-(f_i)\big|_{\spb_r^k}
    = \frac{(-k+m-r)(2+m+r)}{2(r+1)}\,. \qed
\eea
\end{Proposition}

\section{Lower bound of the spectrum of $\Da$ on $P$}

Let $\psi\in \G(\spb(HP))$ be an eigenspinor of $\Da$ with eigenvalue $\l$.
In $p\in P$ the sum over the squares of the 
absolute values of $\cD^{++}_{e_j}\psi$ 
is calculated:
\bea
  \sum_{j=1}^{4m}\|\cD^{++}_{e_j}\psi\|_p^2 &=& 
   2\sum_{j=1}^{2m}\| \Big\la 
    \nba_{f_j} \psi - \sum_{k=1}^{2m}\frac{1}{A_{r+ 1,k+1}^{++}}
    p_{r+ 1}^\mp(f_j)p_r^+(\bar f_k) \nba_{f_k}\psi, \nba_{f_j} \psi \Big\ra
   \nonumber\\ &=&
   2\sum_{j=1}^{2m}\|\nba_{f_j} \psi\|^2 + \frac{1}{2A_{r+ 1,k+1}^{++}}
    \|\Da^{++}\psi\|^2
  \nonumber\\ &\geq & 0\,.
\eea
After division by $2$ and integration over $P$ the
inequality simplifies to
\be
  0 \leq \sum_j\|\nba_{f_j}\psi\|^2 
    + \frac{1}{4 A_{r+1,k+1}^{++}} \|\Da^{++}\psi\|^2\,.
\ee
Analogous inequalities are obtained by dealing with 
$\|\cD_{e_j}^{+-}\psi\|_p^2$, $\|\cD_{e_j}^{-+}\psi\|_p^2$ 
and $\|\cD_{e_j}^{--}\psi\|_p^2$:
\bea
  0 &\leq& \sum_j\|\nba_{\bar f_j}\psi\|^2 
    + \frac{1}{4 A_{r+1,k-1}^{+-}} \|\Da^{+-}\psi\|^2\nonumber\\
  0 &\leq& \sum_j\|\nba_{f_j}\psi\|^2 
    + \frac{1}{4 A_{r-1,k+1}^{-+}} \|\Da^{-+}\psi\|^2\nonumber\\
  0 &\leq& \sum_j\|\nba_{\bar f_j}\psi\|^2 
    + \frac{1}{4 A_{r-1,k-1}^{--}} \|\Da^{--}\psi\|^2\,.
\eea

Now the first and the second resp. the third and the fourth 
inequality are added
and the Weitzenb\"ock formula can be applied:
\bea \label{unglpppm}
  && \sum_j \Big(\|\nba_{f_j}\psi\|^2 + \|\nba_{\bar f_j}\psi\|^2 \Big)
    + \frac{1}{4 A_{r+1,k+1}^{++}} \|\Da^{++}\psi\|^2 
    + \frac{1}{4 A_{r+1,k-1}^{+-}} \|\Da^{+-}\psi\|^2  \nonumber\\
    &&= \frac{1}{2} \Big( \|D\psi\|^2 - \frac{\k}{4}\|\psi\|^2 \Big)
    + \frac{1}{4 A_{r+1,k+1}^{++}} \|\Da^{++}\psi\|^2 
    + \frac{1}{4 A_{r+1,k-1}^{+-}} \|\Da^{+-}\psi\|^2 \nonumber\\
    &&\geq 0
\eea
and in the same manner:
\bea \label{unglmpmm}
  && \sum_j \Big(\|\nba_{f_j}\psi\|^2 + \|\nba_{\bar f_j}\psi\|^2 \Big)
    + \frac{1}{4 A_{r-1,k+1}^{-+}} \|\Da^{-+}\psi\|^2 
    + \frac{1}{4 A_{r-1,k-1}^{--}} \|\Da^{--}\psi\|^2  \nonumber\\
    &&= \frac{1}{2} \Big( \|D\psi\|^2 - \frac{\k}{4}\|\psi\|^2 \Big)
    + \frac{1}{4 A_{r-1,k+1}^{-+}} \|\Da^{-+}\psi\|^2 
    + \frac{1}{4 A_{r-1,k-1}^{--}} \|\Da^{--}\psi\|^2  \nonumber\\
  &&\geq 0\,.
\eea

At this point the problem arises that the expressions $\|\Da^{\pm\pm}\psi\|^2$ 
cannot be calculated directly. Nevertheless it is possible to determine them
by additional assumptions on the eigenspinor $\psi$.

To start with an eigenspinor $\phi$ of $\Da^2$ with eigenvalue $\l^2$, it can
be assumed, that it is localized in a bundle $\spb_r^k$. An
eigenspinor of $\Da$ is determined from this by $\psi := \l \phi \pm \Da \phi$.
Clearly, one has $\psi= \psi_{r-1}+\psi_r+\psi_{r+1} \in \spb_{r-1}
\oplus \spb_r \oplus \spb_r$ and in addition 
$\psi_{r-1}$ and $\psi_{r+1}$ are themselves
eigenspinors of $\Da$. Therefore 
$\l\psi_{r-1} + \Da \psi_{r-1} \in \G(\spb_{r-1}\oplus\spb_{r})$ resp.
$\l\psi_{r+1} + \Da \psi_{r+1}\in \G(\spb_{r+1}\oplus\spb_{r})$ must be
eigenspinors of $\Da$ with eigenvalue $\l$. Hence it can be assumed that an
eigenspinor of $\Da$ is localized in two neighbouring subbundles:
$\psi \in \spb_r\oplus\spb_{r+1}$. If one considers the splitting of $\spb$
w.r.t. $\Oa_1$ instead of $\Oa$, the same 
argumentation is valid. This means, it is possible to assume that for an 
eigenspinor
$\psi$ the following holds: $\psi \in \spb^k \oplus\spb^{k+1}$.
Combining both results, one has:
\be
  \psi \in \G(\spb_r^k \oplus \spb^{k}_{r+1} \oplus \spb^{k+1}_{r} 
  \oplus \spb^{k+1}_{r+1})\,.
\ee
Taking into account the results of Section \ref{Zerlegung}, it can be
seen, that in this direct sum only two summands exist, because the
eigenvalues for $\Oa_1$ of the partial bundles of $\spb_r$ have a distance
of 4. Hence, either case A holds:
\be
  \psi = \psi_0 + \psi_1 \in \G(\spb_r^k \oplus \spb^{k+1}_{r+1}) 
  \quad \hbox{if}\quad
  {\frac{k+r-m}{2}} \in \N{}_0
\ee
or case B:
\be
  \psi = \psi_0 + \psi_1 \in \G(\spb_r^{k+1} \oplus \spb^{k}_{r+1}) 
  \quad \hbox{if}\quad
  {\frac{k+1+r-m}{2}} \in \N{}_0\,.
\ee
Clearly, it is assumed that all the bundles $\spb_r^k$ appearing in the two
cases do exist, i.e. it should not occur that e.g. ${\frac{k+1+r-m}{2}}$ is
an integer but negative. In this case $\psi$ would lie on the allowed lattice
but stick out of the allowed region. Now a picture can be helpful again:

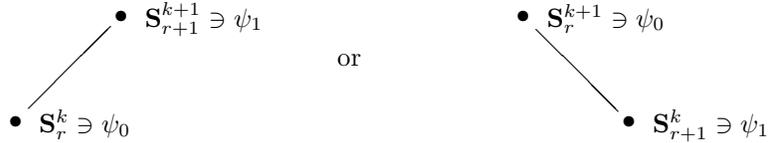
\hbox to \hsize{\hfil
\vbox{\hbox to 100pt{\begin{picture}(100,67)(-10,-10)
 \put(0,0){\circle*{4}}
 \put(9,-4){$\spb_r^k \ni \psi_0$}
 \put(40,40){\circle*{4}}
 \put(49,36){$\spb_{r+1}^{k+1} \ni \psi_1$}
 \put(5,5){\line(1,1){31}}
\end{picture}\hss}}
\vbox to 67pt{\vfil\hbox to 3cm{\hskip1cm or \hss}\vfil}
\vbox{\hbox to 100pt{\begin{picture}(100,67)(-10,-10)
 \put(0,40){\circle*{4}}
 \put(9,36){$\spb_r^{k+1} \ni \psi_0$}
 \put(40,0){\circle*{4}}
 \put(49,-4){$\spb_{r+1}^{k} \ni \psi_1$}
 \put(5,35){\line(1,-1){31}}
\end{picture}\hss}}\hfil}

\begin{Lemma}
In case A one has:
\bea
  \Da^{++}\psi_1=\Da^{+-}\psi_1=\Da^{-+}\psi_1 =0 
   &,& \quad \Da^{--}\psi_1 = \l\psi_0
   \qquad \hbox{and} \nonumber\\
   \Da^{+-}\psi_0=\Da^{-+}\psi_0=\Da^{--}\psi_0=0 
    &,& \quad \Da^{++}\psi_0  = \l \psi_1\,.
\eea
In case B one has:
\bea
  \Da^{++}\psi_1=\Da^{+-}\psi_1=\Da^{--}\psi_1=0 
    &,& \quad \Da^{-+}\psi_1 = \l \psi_0 
   \qquad \hbox{and} \nonumber\\
   \Da^{++}\psi_0=\Da^{-+}\psi_0=\Da^{--}\psi_0=0 
    &,& \quad \Da^{+-}\psi_0 = \l\psi_1\,.
\eea
\end{Lemma}
\proof
This is trivial, because $\psi$ is an eigenspinor of $\Da$, and if one of
the mentioned terms would not vanish, it would be orthogonal to $\psi$.
\leer
With this considerations it is now possible to calculate the estimate of the
eigenvalue.

Case A: Let 
$\psi= \psi_0 + \psi_1 \in \G(\spb_r^k \oplus \spb^{k+1}_{r+1})$ and
${\frac{k+r-m}{2}} \in \N{}_0$. Applying inequality (\ref{unglpppm}) 
to $\psi_0$ leads to (since $\Da^{+-}\psi_0 = 0$):
\be
  \l^2 \geq \frac{2A^{++}_{r+1,k+1}}{2A^{++}_{r+1,k+1}+1}\,\frac{\k}{4}
  = {{\left( -2 - k + m - r \right) \,\left( 3 + m + r \right) }\over 
   {\left( 2 + r \right)  + 
       \left( -2 - k + m - r \right) \,\left( 3 + m + r \right) 
           }} \, \frac{\k}{4}\,.
\ee
In the same manner applying of (\ref{unglmpmm}) to $\psi_1$ leads to
(since $\Da^{-+}\psi_1 = 0$):
\be
  \l^2 \geq \frac{2A^{--}_{r,k}}{2A^{--}_{r,k}+1}\,\frac{\k}{4}
  ={{\left( - m + r \right) \,\left( 2 + k - m + r \right) }\over 
   {r+1 + \left(  - m + r \right) \,
           \left(2 + k - m + r \right)  }}\, \frac{\k}{4}\,.
\ee
Both inequalities must hold simultaneously. It is easy to see that in both
inequalities the right-hand sides are monotonely decreasing with $k$. The
smallest allowed value of $k$ is $m-r$. For this value the expressions are
simplified to:
\be \label{abschFA1}
  \l^2 \geq \frac{2(3+m+r)}{4+2m+r}\,\frac{\k}{4}
\ee
and
\be \label{abschFA2}
  \l^2 \geq \frac{2m-2r}{2m-3r-1}\,\frac{\k}{4} \,.
\ee
The second inequality is weaker than the first, so it can be omitted.
Case B: Let 
$\psi= \psi_0 + \psi_1 \in \G(\spb_r^k \oplus \spb^{k-1}_{r+1})$ and
${\frac{k+r-m}{2}} \in \N{}_0$. After inserting $\psi_0$ into 
(\ref{unglpppm}) one gets (now $\Da^{++}\psi_0 = 0$):
\be
  \l^2 \geq \frac{2A^{+-}_{r+1,k-1}}{2A^{+-}_{r+1,k-1}+1}\,\frac{\k}{4}
  ={{\left( -2 + k - m - r \right) \,\left( 3 + m + r \right) }\over 
   {\left( 2 + r \right)  + 
       \left( -2 + k - m - r \right) \,\left( 3 + m + r \right) 
          }} \frac{\k}{4}\,.
\ee
And finally one has to apply (\ref{unglmpmm}) to $\psi_1$
($\Da^{--}\psi_1 = 0$):
\be
  \l^2 \geq \frac{2A^{--}_{r,k}}{2A^{--}_{r,k}+1}\,\frac{\k}{4}
  = {{\left(  - m + r \right) \,\left( 2-k + m + r \right) }\over 
   {r+1 + \left( - m + r \right) \,
           \left(2 -k + m + r \right)  }}\frac{\k}{4}\,.
\ee
In both inequalities the right-hand sides are monotonely increasing with 
$k$. If the maximal allowed $k=m+r$ is considered, the same inequalities
as (\ref{abschFA1}) and (\ref{abschFA2}) are obtained (if it is taken into
account
that $-\Oa_1$ instead of $\Oa_1$ could have been chosen as distinguished
K\"ahler form, it is clear that case A would have become case B and vice
versa, so both cases have to be equivalent).

\section{The first eigenvalue of $D$ on $M$}

In the preceeding section an estimate of the first eigenvalue of $\Da^2$
on $P$ has been attained. But one is interested in the operator $D^2$ on $M$.
By Proposition \ref{wennDdannDa} it is assured that an eigenvalue of $D$
is an eigenvalue of $\Da$ too, but surely the converse does not hold. 
This means that the estimate for $\Da$ cannot be sharp for $D$.
In addition, the existence of global K\"ahler structures on $P$ has been
used, which do not at all exist on $M$. The right approach is to look only at
horizontal spinors on $P$ that are pull-backs from $M$, i.e. that are of the
form $\pi^\ast \psi$ with $\psi\in \G(\spb(M))$.

If there is given an eigenspinor $\psi \in \G(\spb_r(M))$ one has to 
consider $\tilde\psi := \pi^\ast\psi \in \G(\spb_r)$ first and
to split $\tilde\psi$ into parts which lie in 
$\spb^k_r$ for ${(\frac{k+r-m}{2})} \in \N{}_0$ in order to attain exact
estimates. But it will turn out that $\tilde\psi$ always has contributions
in the subbundles for maximal resp. minimal $k$. Hence the estimate
calculated in the section above for maximal resp. minimal $k$ does hold.

Now let $\psi\in \G(\spb_r)$, $0\leq r\leq m$, be an eigenspinor of $D^2$
on $M$. On $P$ there is a distinguished K\"ahler form $\Oa_1$ operating by
Clifford multiplication on $\tilde\psi$. By definition, in a point 
$p_0\in P$ it holds:
\be
  \Oa_1\tilde\psi(p_0) = \frac{1}{2}\sum_{j=1}^{4m} e^\ast_jJ_1 e^\ast_j 
   \tilde\psi(p_0)
  = \pi^\ast(\Omega_1(p_0) \psi)(p_0) \,,
\ee
where $\Omega_1(p_0):=\frac{1}{2}\sum_{j=1}^{4m} e_jJ_1(p_0) e_j$.
Here one can restrict oneself completely to representation theory. For this,
one considers the splitting of $(\spb_r(M))_{\pi(p_0)}$ w.r.t.
Cliffordmultiplication with $\Omega_1(p_0)$. As seen in Section
\ref{Zerlegung}, $(\spb_r(M))_{\pi(p_0)}
\simeq (\spb_r)_{p_0} = \oplus_{s=0}^r (\spb_r^{(s)})_{p_0}$ splits into 
$r+1$ eigenspaces. There is a corresponding splitting of $\psi$:
\be \label{phizerl}
  \psi= \sum_{s=0}^r \psi_{p_0}^{(s)}\,.
\ee
The index $p_0$ denotes that the splitting depends on the choice of the
three almost-complex structures, i.e. on the point $p_0$ in the fibre of $P$.

Under the operation of $\Oa_1$ the eigenspinor $\tilde\psi$ of $\Da^2$
splits into parts lying in the bundles $\spb_r^{(s)}$. The splitting 
(\ref{phizerl}) corresponds exactly to the splitting of $\tilde\psi$ in
$p_0$. If in (\ref{phizerl}) the summand $\psi_{p_0}^{(0)}$ does not
vanish, it is clear that inequality (\ref{abschFA1}) for the eigenvalue of
$\tilde\psi$ must hold, because on $P$, the subbundles $\spb_r^{(s)}$ 
are invariant under
$\Da^2$. The Dirac eigenvalues of $\tilde\psi$ and of $\psi$ are the same, 
so the
desired estimate holds if $\psi_{p_0}^{(0)}\neq 0$ has been proven.

Now let $p \neq p_0$ be another point on $P$ in the same fibre as $p_0$.
Hence, there exist $g \in \SO(3)$ with $p= p_0 g$. There are some
consequences for the operation of $\Oa_1$: first, as above
\be
  \Oa_1\tilde\psi(p_0g) = \pi^\ast\Big(
  \frac{1}{2}\sum_{j=1}^{4m} e_jJ_1(p_0g) e_j \psi\Big)(p_0g)\,.
\ee
By definition, 
$J_1(p_0g)= g^{-1}J_1(p_0)$ and $\Omega_1(p_0g) = g^{-1}\Omega_1(p_0)$
hold, where $\SO(3)$ operates on $\Omega_1$, $\Omega_2$ and $\Omega_3$ in
the natural way. As above, there is a splitting of $\psi$ corresponding
to the weight space decomposition w.r.t. $g^{-1}\Omega_1(p_0)$:
\be
  \psi= \sum_{s=0}^r \psi_{p_0g}^{(s)}\,.
\ee

In terms of representation theory this means that if $\{H_1,H_2,H_3\}$ forms
the standard base of $\so(3)$ with $\lb H_a,H_b\rb = 2 \e_{abc}H_c$, it is
clear that $\Omega_a$ is the image of $2i H_a$ for $a=1,2,3$ under the
representation determined by the choice of $p_0$ in the given fibre.
Complexification of $\so(3)$ and the definition of 
$X_1 =\frac{1}{2}(H_2 +i H_3)$ and $ Y_1=\frac{1}{2}(H_2 -iH_3)$
leads to the Lie algebra $\sl(2,\C{})$ with the canonical commutator equations.
Hence, $\spb_r(M)_{\pi(p_0)}$ becomes a representation space of $\sl(2,\C{})$,
which will now be identified with an abstract representation space $V$ of
highest weight $r$. This weight is to be regarded as highest eigenvalue of
$H_1$. If an other point $p_0g$ on the fibre is chosen, $r$ is to be regarded
as highest eigenvalue of $g^{-1}H_1$.

For simplicity it can be assumed that $V$ is irreducible. If not, the following
is to be carried out for all irreducible parts separately.

\begin{Proposition}
For every $v\in V$ there is a $g \in \SO(3)$ such that $v$ has contributions
to the subspace of highest weight w.r.t. $gH_1$.
\end{Proposition}
\proof
Let $v$ be given and $g\in \SO(3)$. Let $V$ be equipped with a norm, and let
$\{v^s_g\}$ be a base of normalized vectors spanning the weight spaces
$V^s_g$ w.r.t. $gH_1$. Then $v$ splits into contributions to 
weight spaces w.r.t. $gH_1$ in the following manner:
\be \label{gl1}
  v = \sum_{s=0}^r a^s_g v^s_g, \qquad a^s_g \in \C{}\,.
\ee
Now let $g_t : \lb 0,1\rb \to \SO(3)$ be a path in $\SO(3)$ and $g_0= e$.
By taking derivatives of (\ref{gl1}) one obtains:
\be \label{gl2}
  0 =  \frac{d}{dt}\Big|_{t=0} v =
  \sum_{s=0}^r \Big(\frac{d}{dt}\Big|_{t=0}(a_{g_t}^s)v^s_e + 
   a^s_e \frac{d}{dt}\Big|_{t=0} v_{g_t}^s \Big)\,. 
\ee
In order to calculate the derivative of $v_{g_t}^s$ one considers
\be
  g_tH_1 v_{g_t}^s = (r-2s)v_{g_t}^s
\ee
and derivatives are taken:
\be \label{gl3}
  \frac{d}{dt}\Big|_{t=0}(g_tH_1) v^s_e + 
  H_1 \frac{d}{dt}\Big|_{t=0}(v_{g_t}^s)
  = (r-2s)\frac{d}{dt}\Big|_{t=0} v_{g_t}^s\,.
\ee
Since $g\in\SO(3)$, $\frac{d}{dt}|_{t=0}(g_tH_1) = A_2 H_2 + A_3 H_3$ for
some $A_2,A_3 \in \R{}$. The path $g_t$ can be chosen such that 
$A_2,A_3 \neq 0$.
On the other side, $H_2$ and $H_3$ can be
expressed by the ladder operators $X_1$ and $Y_1$ w.r.t. $H_1$:
\be
  H_2 = X_1 +Y_1 , \qquad H_3 = -i (X_1 -Y_1)\,,
\ee
such that after using (\ref{gl3}) the following holds:
\be
  \big((A_2 -iA_3)X_1 + (A_2 +iA_3)Y_1\big) v^s_e = 
  (r-2s - H_1)\frac{d}{dt}\Big|_{t=0} v_{g_t}^s\,.
\ee
The derivative of $v_{g_t}^s$ in $t=0$ has contributions in $V_e^{s-1}$
and $V_e^{s+1}$ only, explicitly:
\bea
  \Big\lb\frac{d}{dt}\Big|_{t=0} v_{g_t}^s\Big\rb_e^{s-1} 
   &=& \frac{1}{2}(A_2 -iA_3)X_1 v^s_e \nonumber\\
  \Big\lb\frac{d}{dt}\Big|_{t=0} v_{g_t}^s\Big\rb_e^{s+1} 
   &=& -\frac{1}{2}(A_2 +iA_3)Y_1 v^s_e\,.
\eea
The proposition follows easily by contradiction. First it is assumed
that $v$ does not have contributions to the space of highest weight with
respect to $gH_1$ for all $g \in \SO(3)$. Let $s_0>0$ be the largest index
satisfying $a_g^s=0$ for $s<s_0$ and all $g\in \SO(3)$.
Without loss of generality $a^{s_0}_e \neq 0$ can be assumed. But looking
at the contribution to $V_e^{s_0-1}$ in equation (\ref{gl2}):
\be
  0= \frac{d}{dt}\Big|_{t=0}(a_{g_t}^{s_0-1})v_e^{s_0-1}
   +  a_e^{s_0}\cdot\frac{1}{2}(A_2 -iA_3)X_1 v_e^{s_0}\,,
\ee
a contradiction is obtained because by assumption $a_e^{s_0}\neq 0$,
but $a_g^{s_0-1} = 0$ for all $g\in \SO(3)$.
\qed
\leer

It is obvious how to proceed further: the eigenspinor 
$\tilde\psi= \pi^\ast \psi$ has contributions to 
$\spb_r^{(0)}=\spb_r^{m-r}$ and hence
for $\l^2$ in the subbundle $\spb_r(M)$ it holds:
\be
  \l^2 \geq \frac{2(3+m+r)}{4+2m+r}\,\frac{\k}{4}\,.
\ee
The right-hand side is monotonly increasing with $r$, so the universal
estimate is obtained by setting $r=0$:
\be
  \l^2 \geq \frac{m+3}{m+2}\,\frac{\k}{4}\,.
\ee

\section{What comes next?}

There is still the open question for which quaternionic K\"ahler manifolds
of positive scalar curvature the obtained lower bound is sharp.
Maybe the formulation of a hyperk\"ahlerian Killing equation is necessary
and the search for criterions of existence of hyperk\"ahlerian Killing
spinors will give an answer to that question. By analogous methods results
are obtained by C.~B\"ar \cite{Baer4} in the case of Riemannian manifolds and 
by A.~Moroianu \cite{Moroianu1} in the K\"ahlerian case.

\end{document}